\documentclass[11pt]{article}

\usepackage{epsfig,epsf,url,fancybox}
\usepackage{amsmath}
\usepackage{mathrsfs}
\usepackage{amssymb}
\usepackage{graphicx}
\usepackage{color}
\usepackage{changepage}
\usepackage[left=1in,top=1in,right=1in,bottom=1in,letterpaper]{geometry}

\newtheorem{theorem}{Theorem}[section]

\newtheorem{proposition}[theorem]{Proposition}

\newtheorem{example}[theorem]{Example}
\newtheorem{remark}[theorem]{Remark}

\newtheorem{assumption}[theorem]{Assumption}

\newcommand{\Br}{\mathbb{R}}

\newcommand{\st}{\textnormal{s.t.}}

\,
\,

\newcommand{\argmin}{\mathop{\rm argmin}}
\newcommand{\LCal}{\mathcal{L}}

\newcommand{\half}{\frac{1}{2}}

\newcommand{\br}{\mathbb{R}}

\newcommand{\be}{\begin{equation}}
\newcommand{\ee}{\end{equation}}
\newcommand{\ba}{\begin{array}}
\newcommand{\ea}{\end{array}}
\newcommand{\bpm}{\begin{pmatrix}}
\newcommand{\epm}{\end{pmatrix}}
\newcommand{\ACal}{\mathcal{A}}

\newcommand{\XCal}{\mathcal{X}}

\newcommand{\etal}{{et al.\ }}

\begin{document}

\title{Global Convergence of Unmodified 3-Block ADMM for a Class of Convex Minimization Problems}
\author{
Tianyi Lin
\thanks{Department of Industrial Engineering and Operations Research, UC Berkeley, Berkeley, CA 94704. Email: darren\_lin@berkeley.edu.}
\and Shiqian Ma
\thanks{Department of Mathematics, UC Davis, Davis, CA 95616. Email: sqma@math.ucdavis.edu. Research of this author was supported in part by a startup package in Department of Mathematics at UC Davis.}
\and Shuzhong Zhang
\thanks{Department of Industrial and Systems Engineering, University of Minnesota, Minneapolis, MN 55455. Email: zhangs@umn.edu. Research of this author was supported in part by the NSF Grant CMMI-1462408.}}

\date{Nov 7, 2017}

\maketitle

\begin{abstract}

The alternating direction method of multipliers (ADMM) has been successfully applied to solve structured convex optimization problems due to its superior practical performance. The convergence properties of the 2-block ADMM have been studied extensively in the literature. Specifically, it has been proven that the 2-block ADMM globally converges for any penalty parameter $\gamma>0$. In this sense, the 2-block ADMM allows the parameter to be free, i.e., there is no need to restrict the value for the parameter when implementing this algorithm in order to ensure convergence. However, for the 3-block ADMM, Chen \etal \cite{Chen-admm-failure-2013} recently constructed a counter-example showing that it can diverge if no further condition is imposed. The existing results on studying further sufficient conditions on guaranteeing the convergence of the 3-block ADMM usually require $\gamma$ to be smaller than a certain bound, which is usually either difficult to compute or too small to make it a practical algorithm. In this paper, we show that the 3-block ADMM still globally converges with any penalty parameter $\gamma>0$ if the third function $f_3$ in the objective is smooth and strongly convex, and its condition number is in $[1,1.0798)$, besides some other mild conditions. This requirement covers an important class of problems to be called regularized least squares decomposition (RLSD) in this paper.

\vspace{0.5cm}




\end{abstract}

\section{Introduction}

The alternating direction method of multipliers (ADMM) has been very successfully applied to solve many structured convex optimization problems arising from machine learning, image processing, statistics, computer vision and so on; see the recent survey paper \cite{Boyd-etal-ADM-survey-2011}. The ADMM is particularly efficient when the problem has a separable structure in functions and variables. For example, the following convex minimization problem with 2-block variables can usually be solved by ADMM, provided that a certain structure of the problem is in place:
\be\label{sum-2}\ba{ll} \min & f_1(x_1) + f_2(x_2) \\ \st & A_1x_1 + A_2x_2 = b \\ & x_1\in\XCal_1, x_2\in\XCal_2,\ea\ee
where $f_{i}(x_{i}):\Br^{n_{i}}\rightarrow\Br, i=1,2$, are closed convex functions, $A_{i}\in\Br^{p\times n_{i}}, i=1,2$, $b\in\Br^{p}$ and $\XCal_i, i=1,2$, are closed convex sets. A typical iteration of the 2-block ADMM (with given $(x_2^k,\lambda^k)$) for solving \eqref{sum-2}
can be described as
\be\label{admm-2}\left\{\ba{ll} x_1^{k+1} & := \argmin_{x_1\in\XCal_1} \ \bar\LCal_\gamma(x_1,x_2^k;\lambda^k) \\
                                x_2^{k+1} & := \argmin_{x_2\in\XCal_2} \ \bar\LCal_\gamma(x_1^{k+1},x_2;\lambda^k) \\
                                \lambda^{k+1} & := \lambda^k - \gamma(A_1x_1^{k+1}+A_2x_2^{k+1}-b), \ea\right.\ee
where the augmented Lagrangian function $\bar\LCal_\gamma$ is defined as
\[\bar\LCal_\gamma(x_1,x_2;\lambda) := f_1(x_1) + f_2(x_2)-\langle\lambda, A_1x_1 + A_2x_2 - b\rangle + \frac{\gamma}{2}\|A_1x_1 + A_2x_2 - b\|_2^2,\]
where $\lambda$ is the Lagrange multiplier and $\gamma>0$ is a penalty parameter, which can also be viewed as a step size on the dual update. The convergence properties of 2-block ADMM \eqref{admm-2} have been studied extensively in the literature; see for example \cite{Lions-Mercier-79,Gabay-83,Fortin-Glowinski-1983,Glowinski-LeTallec-89,Eckstein-Bertsekas-1992,He-Yuan-rate-ADM-2012,Monteiro-Svaiter-2010a,Deng-Yin-2012,Boley-2012}. A very nice property of the 2-block ADMM is that it is {\it parameter restriction-free}: it has been proven that the 2-block ADMM \eqref{admm-2} is globally convergent for any parameter $\gamma>0$, starting from anywhere. This property makes the 2-block ADMM particularly attractive for solving structured convex optimization problems in the form of \eqref{sum-2}.

However, this is not the case when ADMM is applied to solve convex problems with 3-block variables:
\be\label{sum-3}\ba{ll} \min & f_1(x_1) + f_2(x_2) + f_3(x_3) \\ \st & A_1x_1 + A_2x_2 + A_3x_3 = b \\ & x_1\in\XCal_1, x_2\in\XCal_2,x_3\in\XCal_3.\ea\ee
Note that the 3-block ADMM for solving \eqref{sum-3} can be described as
\be\label{admm-3}\left\{\ba{ll} x_1^{k+1} & := \argmin_{x_1\in\XCal_1} \ \LCal_\gamma(x_1,x_2^k,x_3^k;\lambda^k) \\
                                x_2^{k+1} & := \argmin_{x_2\in\XCal_2} \ \LCal_\gamma(x_1^{k+1},x_2,x_3^k;\lambda^k) \\
                                x_3^{k+1} & := \argmin_{x_3\in\XCal_3} \ \LCal_\gamma(x_1^{k+1},x_2^{k+1},x_3;\lambda^k) \\
                                \lambda^{k+1} & := \lambda^k - \gamma(A_1x_1^{k+1}+A_2x_2^{k+1}+A_3x_3^{k+1}-b), \ea\right.\ee
where the augmented Lagrangian function is defined as
\[\LCal_\gamma(x_1,x_2,x_3;\lambda) := f_1(x_1) + f_2(x_2) + f_3(x_3)-\langle\lambda, A_1x_1 + A_2x_2+A_3x_3 - b\rangle + \frac{\gamma}{2}\|A_1x_1 + A_2x_2+A_3x_3 - b\|_2^2.\]
Regarding its general convergence however, Chen \etal constructed a counterexample in \cite{Chen-admm-failure-2013} showing that the 3-block ADMM \eqref{admm-3} can diverge if no further condition is imposed. On the other hand, the 3-block ADMM \eqref{admm-3} has been successfully used in many important applications such as the robust and stable principal component pursuit problem \cite{Tao-Yuan-SPCP-2011,Zhou-Candes-Ma-StablePCP-2010}, the robust image alignment problem \cite{Wright-RASL-TPAMI}, Semidefinite Programming \cite{Wen-Goldfarb-Yin-2009}, and so on. It is therefore of great interest to further study sufficient conditions to guarantee the convergence of 3-block ADMM \eqref{admm-3}. Han and Yuan \cite{Han-Yuan-note-2012} showed that the 3-block ADMM \eqref{admm-3} converges if all the functions $f_1,f_2,f_3$ are strongly convex and $\gamma$ is restricted to be smaller than a certain bound. This condition is relaxed in Chen, Shen and You \cite{Chen-Shen-You-convergence-2013} and Lin, Ma and Zhang \cite{Lin-Ma-Zhang-convergence-2014} to allow only $f_2$ and $f_3$ to be strongly convex and $\gamma$ is restricted to be smaller than a certain bound. Moreover, the first sublinear convergence rate result of multi-block ADMM is established in \cite{Lin-Ma-Zhang-convergence-2014}. Closely related to \cite{Chen-Shen-You-convergence-2013, Lin-Ma-Zhang-convergence-2014}, Cai, Han and Yuan \cite{Cai-Han-Yuan-direct-2014} and Li, Sun and Toh \cite{Li-Sun-Toh-convergent-2015} proved the convergence of the 3-block ADMM \eqref{admm-3} under the assumption that only one of the functions $f_1$, $f_2$ and $f_3$ is strongly convex, and $\gamma$ is restricted to be smaller than a certain bound. Davis and Yin \cite{Davis-Yin-2015} studied a variant of the 3-block ADMM (see Algorithm 8 in \cite{Davis-Yin-2015}) which requires that $f_1$ is strongly convex and $\gamma$ is smaller than a certain bound to guarantee the convergence. In addition to strong convexity of $f_2$ and $f_3$, and the boundedness of $\gamma$, by assuming further conditions on the smoothness of the functions and some rank conditions on the matrices in the linear constraints, Lin, Ma and Zhang \cite{Lin-Ma-Zhang-2014-linear} proved the globally linear convergence of 3-block ADMM \eqref{admm-3}. More recently, Lin, Ma and Zhang \cite{Lin-Ma-Zhang-2015} further proposed several alternative approaches to ensure the sublinear convergence rate of \eqref{admm-3} without requiring any function to be strongly convex.  Remark that in all these works, to trade for a convergence guarantee the penalty parameter $\gamma$ is required to be small, which potentially affects the practical effectiveness of the 3-block ADMM \eqref{admm-3}, while the 2-block ADMM \eqref{admm-2} 
does not suffer from such compromises.

Alternatively, one may opt to modify the 3-block ADMM \eqref{admm-3} to achieve convergence,
with similar per-iteration computational complexity as \eqref{admm-3}.
The existing methods in the literature along this line can be classified into the following three main categories. (i) The first class of algorithms requires a correction step in the updates (see, e.g., \cite{He-Tao-Yuan-IMA-2014,He-Tao-Yuan-2012,He-Tao-Yuan-MOR-2013,He-Hou-Yuan-Jacob-2013}).
(ii) The second class of algorithms adds proximal terms and/or dual step size to the ADMM updates, i.e., these algorithms change \eqref{admm-3} to
\be\label{admm-3-proximal}\left\{\ba{ll} x_1^{k+1} & := \argmin_{x_1\in\XCal_1} \ \LCal_\gamma(x_1,x_2^k,x_3^k;\lambda^k) + \half\|x-x_1^k\|_{P_1}\\
                                x_2^{k+1} & := \argmin_{x_2\in\XCal_2} \ \LCal_\gamma(x_1^{k+1},x_2,x_3^k;\lambda^k) + \half\|x-x_2^k\|_{P_2} \\
                                x_3^{k+1} & := \argmin_{x_3\in\XCal_3} \ \LCal_\gamma(x_1^{k+1},x_2^{k+1},x_3;\lambda^k) + \half\|x-x_3^k\|_{P_3} \\
                                \lambda^{k+1} & := \lambda^k - \alpha\gamma(A_1x_1^{k+1}+A_2x_2^{k+1}+A_3x_3^{k+1}-b), \ea\right.\ee
where matrices $P_i\succeq 0$ and $\alpha>0$ denotes a step size for the dual update. Global convergence and convergence rate for \eqref{admm-3-proximal} and its variants (for example, allowing to update $x_1$, $x_2$, $x_3$ in a Jacobian manner instead of a Gauss-Seidel manner) are analyzed under various conditions (see, e.g., \cite{Luo-ADMM-2012,Deng-admm-2014,Hong-etal-2014-BSUMM,Sun-Toh-Yang-admm4-2014,Li-Sun-Toh-convergent-2015}).
Note that these works usually require restrictive conditions on $P_i$, $\alpha$ and $\gamma$ that may also affect the performance of solving large-scale problems arising from practice.
Notwithstanding all these efforts, many authors acknowledge
that the unmodified 3-block ADMM \eqref{admm-3} usually outperforms its variants \eqref{admm-3-proximal} and the ones with correction step in practice (see, e.g., the discussions in \cite{Sun-Toh-Yang-admm4-2014,Wang-etal-2013-multi-2}). (iii) The recent work by Sun, Luo and Ye \cite{Sun-Luo-Ye-random-admm-2015} on a randomly permuted ADMM is probably the only variant of 3-block ADMM which does not restrict the $\gamma$ value, but its convergence
is now only guaranteed for solving a squared and nonsingular linear system.

Motivated by the fact that the 2-block ADMM \eqref{admm-2} allows the parameter to be free,
in this paper we set out to explore the structures of 3-block model for which the unmodified 3-block ADMM \eqref{admm-3} converges for all parameter values. Given the superior performance of \eqref{admm-3},
such property is of great practical importance.
In this paper, we show that the 3-block ADMM \eqref{admm-3} is globally convergent for any fixed $\gamma>0$ if $A_3=I$, $x_3\in\mathcal{X}_3$ is absent, and more importantly, the condition number of $f_3$ is in $[1,1.0798)$, which covers an important class of convex problems, termed the Regularized Least Squares Decomposition (RLSD) in this paper.

\section{Preliminaries}

In this paper, we consider the 3-block ADMM for solving
\be\label{prob:P-uncons}
\ba{ll}
\min & f_1(x_1) + f_2(x_2) + f_3(x_3) \\
\st & A_1 x_1 + A_2 x_2 + x_3 = b, \quad x_i\in\XCal_i, \ i=1,2.
\ea
\ee
It is noted that \eqref{prob:P-uncons} is a special case of \eqref{sum-3} with $A_3$ being identity and $x_3\in\mathcal{X}_3$ being removed. Throughout this paper, we make the following assumption on $f_3$.
\begin{assumption} \label{assumption_3}
We assume that function $f_3$ is lower bounded by $f_3^*$ and is strongly convex with modulus $\sigma>0$ and $\nabla f_3$ is Lipschitz continuous with Lipschitz constant $L>0$; i.e., the following inequalities hold:
\[\inf_{x_3\in\br^p} f_3(x_3) > f_3^* > -\infty, \]
\be\label{strong_convex_1}
f_3(y) \geq f_3(x) + \left( y-x \right)^{\top}\nabla f_3(x) + \frac{\sigma}{2} \left\|y - x \right\|^2, \quad \forall x,y\in\br^p
\ee
or equivalently,
\be\label{strong_convex_2}
(y-x)^{\top}\left( \nabla f_3(y) - \nabla f_3(x) \right) \geq \sigma\left\| y - x \right\|^2, \quad \forall x,y\in\br^p
\ee
and
\be\label{lipschitz_continuous}
\left\| \nabla f_3(y) - \nabla f_3(x) \right\| \leq L \left\| y - x \right\|, \quad \forall x,y\in\br^p.
\ee
\end{assumption}

A special case of $f_3$ that satisfies Assumption \ref{assumption_3} is $f_3(\cdot)=\half\|\cdot\|_2^2$. In this case, \eqref{prob:P-uncons} can be equivalently written as:
\be\label{sharing}\ba{ll} \min & f_1(x_1) + f_2(x_2) + \half\left\|A_1 x_1 + A_2 x_2-b\right\|^2  \\
\st &  x_1\in\XCal_1, x_2\in\XCal_2.\ea\ee
We call \eqref{sharing} regularized least squares decomposition (RLSD).
In \eqref{sharing}, one seeks to decompose the observed data $b$ into two components $A_1x_1$ and $A_2x_2$, and $f_1$ and $f_2$ denote some regularization functions that promote certain structures of $x_1$ and $x_2$ in the decomposed terms.
One may also view \eqref{sharing} as a data fitting problem with two regularization terms,
where $\left\|A_1 x_1 + A_2 x_2-b\right\|^2$ denotes a least squares loss function on the data fitting term.
Many works in the literature (including Boyd \etal \cite{Boyd-etal-ADM-survey-2011} and Hong, Luo and Razaviyayn \cite{Hong-nonconvex-admm-2014}) have suggested to solve \eqref{sharing} by applying ADMM to solve its equivalent form \eqref{prob:P-uncons}. The advantage of using ADMM to solve \eqref{prob:P-uncons} is that the subproblems are usually easy to solve. Especially, the subproblem for $x_3$ has a closed-form solution. Yang and Zhang \cite{Yang-Zhang-2009} applied the 2-block ADMM to solve the following $\ell_1$-norm regularized least squares problem (or the so-called Lasso problem \cite{Tibshirani-LASSO-1996} in statistics):
\be\label{lasso} \min_x \ \beta\|x\|_1 + \half\|Ax-b\|^2, \ee
where $\beta>0$ is a weighting parameter. Therefore, the Lasso problem is in fact RLSD with one block of variables (more on this later). In order to use ADMM, Yang and Zhang \cite{Yang-Zhang-2009} reformulated \eqref{lasso} as
\be\label{lasso-reformulate}\ba{ll} \min_{x,r} & \beta\|x\|_1 + \half\|r\|^2 \\ \st & Ax-r=b,\ea\ee
in which the two-block variables $x$ and $r$ are associated with two structured functions $\|x\|_1$ and $\|r\|^2$, respectively. Numerical experiments conducted in \cite{Yang-Zhang-2009} showed that the 2-block ADMM greatly outperforms other state-of-the-art solvers on this problem. It is noted that the problem RLSD \eqref{sharing} reduces to the Lasso problem \eqref{lasso} when $f_2$ and $x_2$ vanish and $f_1$ is the $\ell_1$ norm.
Problem RLSD \eqref{sharing} actually covers many interesting applications in practice, and in the following we will discuss a few examples. RLSD \eqref{sharing} is sometimes also known
as sharing problem in the literature, and we refer the interested readers to \cite{Boyd-etal-ADM-survey-2011} and \cite{Hong-nonconvex-admm-2014} for more examples of this problem.

\begin{example}\label{example-1}
Stable principal component pursuit \cite{Zhou-Candes-Ma-StablePCP-2010}. This problem aims to recover a low-rank matrix (the principal components) from a high dimensional data matrix despite both small entry-wise noise and gross sparse errors. This problem can be formulated as (see Eq.\ (15) of \cite{Zhou-Candes-Ma-StablePCP-2010}):
\be\label{spcp}\min_{L,S} \ \beta_1\|L\|_* + \beta_2\|S\|_1 + \half\|M-L-S\|_F^2,\ee
where $M\in\Br^{m\times n}$ is the given corrupted data matrix, $L$ and $S$ are respectively low-rank and sparse component of $M$. It is obvious that this problem is in the form of \eqref{sharing} with $\XCal_1=\XCal_2=\Br^{m\times n}$. For solving \eqref{spcp} using the 3-block ADMM \eqref{admm-3}, see \cite{Tao-Yuan-SPCP-2011}.
\end{example}

\begin{example}\label{example-2}
Static background extraction from surveillance video \cite{Yuan-2015-statis-background-extraction-2015,Ma-SparCoc-2014}. This problem aims to extract the static background from a surveillance video. Given a sequence of frames of a surveillance video $M\in\Br^{m\times n}$, this problem finds a decomposition of $M$ in the form of $M=ue^\top + S$, where $u\in\Br^m$ denotes the static background of the video, $e$ is the all-ones vector, and $S$ denotes the sparse moving foreground in the video. Since the components of u represent the pixel values of the background image, we can restrict $u$ as $b_\ell \leq u \leq b_u$, with $b_\ell=0$ and $b_u=255$. This problem can then be formulated as
\be\label{static-background} \ba{ll}\min_{u, S} & \beta\|S\|_1 + \half\|M-ue^\top - S\|_F^2 \\
                                    \st  & b_\ell \leq u \leq b_u. \ea\ee
Note that \eqref{static-background} is a slight modification of Eq.\ (1.9) in \cite{Yuan-2015-statis-background-extraction-2015} with the bounded constraints added to $u$ in order to get a background image with physical meanings. A similar model was considered by Ma \etal in \cite{Ma-SparCoc-2014} for molecular pattern discovery and cancer gene identification. We refer the interested readers to \cite{Yuan-2015-statis-background-extraction-2015} and \cite{Ma-SparCoc-2014} for more details of this problem.
\end{example}

\begin{example}\label{example-3}
Compressive Principal Component Pursuit \cite{Wright-compressive-pca-2012}. This problem also considers decomposing a matrix $M$ into a low-rank part and a sparse part as \eqref{spcp}. The difference is that $M$ is observed via a small set of linear measurements. This problem can thus be formulated as
\be\label{cpcp}\min_{L,S} \ \beta_1\|L\|_* + \beta_2\|S\|_1 + \half\|M-\ACal(L)-\ACal(S)\|_F^2,\ee
where $\ACal:\Br^{m\times n}\rightarrow \Br^{m\times n}$ is a linear mapping. Note that \eqref{cpcp} is an unconstrained version of Eq.\ (1.7) in \cite{Wright-compressive-pca-2012}, and \eqref{cpcp} is particularly interesting when there are noises in the compressive measurements $M$. Similar problem has also been considered in \cite{WSB-NIPS-2011}.
\end{example}

In this paper, we prove that the unmodified 3-block ADMM \eqref{admm-3} globally converges with any parameter $\gamma>0$, when it is applied to solve problem \eqref{prob:P-uncons}, if $f_3$ satisfies Assumption \ref{assumption_3} and its condition number is in $[1,1.0798)$. This result provides theoretical foundations for using the unmodified 3-block ADMM with {\it a free choice of any parameter $\gamma>0$}.

According to the first-order optimality conditions for \eqref{prob:P-uncons}, solving \eqref{prob:P-uncons} is equivalent to finding $x_1^*\in\XCal_1$, $x_2^*\in\XCal_2$, $x_3^*\in\Br^p$ and $\lambda^*\in\Br^p$
such that the following holds:
\begin{equation}\label{kkt}
\left\{
\begin{array}{l}
f_1(x_1) - f_1(x_1^*) - \left(x_1-x_1^*\right)^\top \left(A_1^\top\lambda^* \right) \geq 0, \quad \forall x_1\in\XCal_1,\\
f_2(x_2) - f_2(x_2^*) - \left(x_2-x_2^*\right)^\top \left(A_2^\top\lambda^* \right) \geq 0, \quad \forall x_2\in\XCal_2,\\
\nabla f_3(x_3^*) - \lambda^* = 0, \\
A_1 x_1^* + A_2 x_2^* + x_3^* = b.
\end{array}
\right.
\end{equation}
We call $(x_1^*,x_2^*,x_3^*)$ optimal primal solution, and $\lambda^*$ optimal dual solution of \eqref{prob:P-uncons}. We use $\Omega^*$ to denote the set of optimal primal and dual pairs $(x_1^*,x_2^*,x_3^*,\lambda^*)$.

The following two assumptions are made throughout this paper.

\begin{assumption}\label{assumption-1}
The set of optimal primal and dual pairs of problem \eqref{prob:P-uncons}, $\Omega^*$, is non-empty.
\end{assumption}

\begin{assumption} \label{assumption-2}
We assume the following conditions hold.
\begin{enumerate}
\item $A_1$ and $A_2$ have full column rank.
\item The objective functions $f_1$ and $f_2$ are lower semi-continuous, and proper closed convex functions.
\item $f_i+\textbf{1}_{\XCal_i}, i=1,2$, are both coercive functions, where $\textbf{1}_{\XCal_i}$ denotes the indicator function of $\XCal_i$, i.e.,
\[\textbf{1}_{\XCal_i}(x_i) = \left\{\ba{ll} 0, & \mbox{ if } x_i\in \XCal_i \\ +\infty, & \mbox{ otherwise. }\ea\right.\]
Note that this assumption implies that $f_1$ and $f_2$ have finite lower bounds on $\XCal_1$ and $\XCal_2$, respectively, i.e.,
\begin{displaymath}
\inf_{x_1\in\XCal_1} f_1(x_1) > f_1^* > -\infty, \quad \inf_{x_2\in\XCal_2} f_2(x_2) > f_2^* > -\infty.
\end{displaymath}
\end{enumerate}
\end{assumption}

\begin{remark}
We remark here that requiring $f_i+\textbf{1}_{\XCal_i}$ to be a coercive function is not a restrictive assumption. Many functions used as regularization terms including $\ell_1$-norm, $\ell_2$-norm, $\ell_\infty$-norm for vectors and nuclear norm for matrices are all coercive functions; assuming the compactness of $\XCal_i$ also leads to the coerciveness of $f_i+\textbf{1}_{\XCal_i}$. For instance, problems considered in Examples \ref{example-1}-\ref{example-3} all satisfy this assumption.
\end{remark}

The following assumption will be used in Theorem \ref{thm} for proving a stronger convergence result.
\begin{assumption}\label{assumption-4}
One of the following two cases holds:
\begin{enumerate}
\item Case (i): $\left[ A_1 \ A_2\right]$ is of full column rank;
\item Case (ii): For $i=1,2$, $f_i$ is locally strongly convex, i.e., there exists $\sigma_i>0$, such that
\begin{equation*}
f_i(x_i) - f_i(x_i^*) - \left( x_i - x_i^*\right)^\top g_i(x_i^*) \geq \frac{\sigma_i}{2}\left\| x_i - x_i^*\right\|^2, \forall x_i\in\XCal_i.
\end{equation*}
\end{enumerate}
\end{assumption}

In our analysis, the following well-known identity and inequality are used frequently:
\begin{eqnarray}
 (w_{1}-w_{2})^\top(w_{3}-w_{1}) &=& \frac{1}{2}\left(\|w_{2}-w_{3}\|^{2}-\|w_{1}-w_{2}\|^{2}-\|w_{1}-w_{3}\|^{2}\right), \label{identity-3} \\
 w_1^\top w_2 &\geq & -\frac{1}{2\xi} \left\| w_1 \right\|^2 - \frac{\xi}{2}\left\| w_2\right\|^2, \quad \forall\xi>0. \label{identity-2}
\end{eqnarray}

\section{Global convergence of 3-block ADMM}

In this section, we show that the 3-block ADMM \eqref{admm-3} converges under Assumptions \ref{assumption_3}, \ref{assumption-1}, \ref{assumption-2} and \ref{assumption-4}, when it is applied to solve \eqref{prob:P-uncons}, given that $\gamma$ is chosen to be any value in the following range:
\begin{eqnarray}\label{assumption}
\gamma & \in & \left(0, \min\left\{\frac{4\sigma}{\eta_2}, \frac{\sigma(\eta_2 - 2)}{4\eta_2} + \sqrt{\frac{\sigma^2(\eta_2 - 2)^2}{16\eta_2^2} + \frac{\sigma^2(\eta_2 - 2)}{4\eta_2}}\right\}\right) \bigcup\left(\sqrt{\sigma^2 + \frac{2L^2}{\eta_1 - 2}} - \sigma, \frac{4\sigma}{\eta_1}\right] \nonumber \\
& & \bigcup\left(\frac{L^2}{\sigma}, +\infty\right),
\end{eqnarray}
where $\eta_1$ and $\eta_2$ can be any value in $(2,+\infty)$. Note that if $\eta_1$ is chosen such that $\sqrt{\sigma^2 + \frac{2L^2}{\eta_1 - 2}} - \sigma > \frac{4\sigma}{\eta_1}$, then the second interval in \eqref{assumption} is empty.

\begin{remark}\label{remark:gamma-bound}
In Proposition \ref{prop:gamma-free} we show that if the condition number of $f_3$, which is defined as $L/\sigma$, is in the interval $[1,1.0798)$, then we can always find $\eta_1>2$ and $\eta_2>2$ such that \eqref{assumption} reduces to $\gamma\in(0,\infty)$. Moreover, the bound $1.0798$ is unlikely to be tight, but we leave this as an open question for future research.
\end{remark}

Note that the 3-block ADMM for solving \eqref{prob:P-uncons} can be written as
\be\label{admm-3-extend}\left\{\ba{ll} x_1^{k+1} & := \argmin_{x_1\in\XCal_1} \ f_1(x_1) + \frac{\gamma}{2}\|A_1x_1+A_2x_2^{k+1}+x_3^{k+1}-b-\lambda^k/\gamma\|^2 \\
                                x_2^{k+1} & := \argmin_{x_2\in\XCal_2} \ f_2(x_2) + \frac{\gamma}{2}\|A_1x_1^{k+1}+A_2x_2+x_3^{k+1}-b-\lambda^k/\gamma\|^2  \\
                                x_3^{k+1} & := \argmin_{x_3\in\br^p} \ f_3(x_3) + \frac{\gamma}{2}\|A_1x_1^{k+1}+A_2x_2^{k+1}+x_3-b-\lambda^k/\gamma\|^2 \\
                                \lambda^{k+1} & := \lambda^k - \gamma(A_1x_1^{k+1}+A_2x_2^{k+1}+x_3^{k+1}-b). \ea\right.\ee
The first-order optimality conditions for the three subproblems in \eqref{admm-3-extend} are given by $x_i^{k+1}\in\XCal_i$ and $x_i\in\XCal_i$ for $i=1,2$, and
\begin{align}
& \left(x_1 - x_1^{k+1}\right)^\top\left[ g_1(x_1^{k+1})-A_1^\top\lambda^k+\gamma A_1^\top \left(A_1 x_1^{k+1} + A_2 x_2^k + x_3^k - b \right)\right] \geq 0, \label{opt-x-1} \\
& \left(x_2 - x_2^{k+1}\right)^\top \left[ g_2(x_2^{k+1})-A_2^\top\lambda^k+\gamma A_2^\top\left(A_1 x_1^{k+1} + A_2 x_2^{k+1} + x_3^k -b\right)\right] \geq 0, \label{opt-x-2} \\
& \nabla f_3(x_3^{k+1})-\lambda^k + \gamma \left(A_1 x_1^{k+1} + A_2 x_2^{k+1} + x_3^{k+1} -b\right) = 0, \label{opt-x-3}
\end{align}
where $g_i \in \partial f_i$ is the subgradient of $f_i$ for $i=1,2$. Moreover, by combining with the updating formula for $\lambda^{k+1}$, \eqref{opt-x-1}-\eqref{opt-x-3} can be rewritten as
\begin{align}
& \left(x_1 - x_1^{k+1}\right)^\top\left[g_1(x_1^{k+1}) - A_1^\top\lambda^{k+1} +\gamma A_1^\top\left( A_2(x_2^k - x_2^{k+1})+ (x_3^k - x_3^{k+1})\right)\right] \geq 0, \label{opt-x-1-lambda} \\
& \left(x_2 - x_2^{k+1}\right)^\top\left[g_2(x_2^{k+1}) - A_2^\top\lambda^{k+1} +\gamma A_2^\top\left( x_3^k - x_3^{k+1}\right)\right] \geq 0, \label{opt-x-2-lambda} \\
& \nabla f_3(x_3^{k+1}) - \lambda^{k+1} = 0. \label{opt-x-3-lambda}
\end{align}

We are now ready to present the main result on the convergence of 3-block ADMM.

\begin{theorem}\label{thm}
Assume Assumptions \ref{assumption-1} and \ref{assumption-2} hold.
Let $\left(x_1^k,x_2^k,x_3^k,\lambda^k\right)$ be generated by the 3-block ADMM \eqref{admm-3-extend} applied to \eqref{prob:P-uncons} with $f_3$ satisfying Assumption \ref{assumption_3} and $\gamma$ chosen as in \eqref{assumption}.
The following results hold.
\begin{enumerate}
\item If
\be\label{gamma-1}
\gamma\in (L^2/\sigma, +\infty),
\ee
then
$\left\{\left(x_1^k,x_2^k,x_3^k,\lambda^k\right): k=0,1,2,\ldots\right\}$ is a bounded sequence and any of its cluster point is an optimal primal and dual pair of \eqref{prob:P-uncons}. Moreover, we have
\begin{eqnarray}\label{thm-obj-feasibility}
\lim_{k\to\infty}\left| f(x_1^k) + f_2(x_2^k) + f_3(x_3^k) - f^* \right| = 0, \quad
\lim_{k\to\infty}\left\| A_1 x_1^k + A_2 x_2^k + x_3^k - b \right\| = 0,
\end{eqnarray}
where $f^*$ denotes the optimal objective value of problem \eqref{prob:P-uncons}. Additionally, if Assumption \ref{assumption-4} holds, then the whole sequence of $\left\{(x_1^k,x_2^k,x_3^k,\lambda^k): k=0,1,2,\ldots\right\}$ converges to an optimal primal and dual pair of problem \eqref{prob:P-uncons}.
\item If
\be\label{gamma-2}
\gamma\in\left(\sqrt{\sigma^2 + \frac{2L^2}{\eta_1 - 2}} - \sigma, \frac{4\sigma}{\eta_1}\right] \bigcup \left(0, \min\left\{\frac{4\sigma}{\eta_2}, \frac{\sigma(\eta_2 - 2)}{4\eta_2} + \sqrt{\frac{\sigma^2(\eta_2 - 2)^2}{16\eta_2^2} + \frac{\sigma^2(\eta_2 - 2)}{4\eta_2}}\right\}\right)
\ee
with $\eta_1$ and $\eta_2$ arbitrarily chosen in $(2,+\infty)$, then
    $\left\{\left(x_1^k,x_2^k,x_3^k,\lambda^k\right): k=0,1,2,\ldots\right\}$ is a bounded sequence,
    and the whole sequence of $\left\{\left(x_1^k,x_2^k,x_3^k,\lambda^k\right): k=0,1,2,\ldots\right\}$ converges to an optimal primal and dual pair of problem \eqref{prob:P-uncons}.
\end{enumerate}
\end{theorem}

\begin{proof}
By \eqref{opt-x-3-lambda} and the Lipschitz continuity of $\nabla f_3$, we have
\begin{equation}\label{lambda-bound-x}
\| \lambda^{k+1} - \lambda^{k} \| \leq L \| x_3^{k+1} - x_3^k \|.
\end{equation}

Letting $x_2 = x_2^k$ in the $(k+1)$-th iteration and $x_2 = x_2^{k+1}$ in the $k$-th iteration of \eqref{opt-x-2-lambda} yields
\begin{eqnarray*}
(x_2^k - x_2^{k+1})^\top\left[ g_2(x_2^{k+1}) - A_2^\top\lambda^{k+1} + \gamma A_2^\top\left( x_3^k - x_3^{k+1}\right)\right] & \geq & 0, \\
(x_2^{k+1}-x_2^k)^\top\left[ g_2(x_2^k) - A_2^\top\lambda^k + \gamma A_2^\top \left(  x_3^{k-1} - x_3^k\right)\right] & \geq & 0.
\end{eqnarray*}
Adding these two inequalities, using the monotonicity of $g_2$ and applying \eqref{identity-2} we obtain that the following inequality holds for any $\epsilon>0$:
\begin{eqnarray}\label{inequality-x-2-primal-dual}
& & \left(A_2 x_2^{k+1} - A_2 x_2^k\right)^{\top}\left( \lambda^{k+1} - \lambda^k \right) \nonumber \\
& \geq & -\frac{\gamma}{\epsilon}\left\| A_2 x_2^{k+1} - A_2 x_2^k \right\|^2 - \frac{\gamma\epsilon}{2}\left\| x_3^k - x_3^{k-1}\right\|^2 - \frac{\gamma\epsilon}{2}\left\| x_3^{k+1} - x_3^k\right\|^2.
\end{eqnarray}
From \eqref{opt-x-3-lambda} and the strong convexity of $f_3$, we have
\begin{eqnarray}\label{inequality-x-3-primal-dual}
\left( x_3^{k+1} - x_3^k\right)^{\top}\left( \lambda^{k+1} - \lambda^k \right) \geq \sigma \left\| x_3^{k+1} - x_3^k \right\|^2 .
\end{eqnarray}

Now we prove part 1. Firstly, we prove that the augmented Lagrangian function $\LCal_\gamma(x_1^k,x_2^k,x_3^k;\lambda^k)$ is non-increasing. Note that the augmented Lagrangian function of \eqref{prob:P-uncons} is
\[
\LCal_\gamma(x_1,x_2,x_3;\lambda) = f_1(x_1) + f_2(x_2) + f_3(x_3) -\langle\lambda,A_1x_1+A_2x_2+x_3-b\rangle+\frac{\gamma}{2}\|A_1x_1+A_2x_2+x_3-b\|_2^2.
\]
Following the same steps as in the proof of Lemma 2.2 of \cite{Hong-nonconvex-admm-2014}, we get the following inequality:
\begin{eqnarray}\label{L-decrease-M}
& & \LCal_\gamma(x_1^k,x_2^k, x_3^k,\lambda^k) - \LCal_\gamma(x_1^{k+1},x_2^{k+1}, x_3^{k+1},\lambda^{k+1}) \nonumber \\
& \geq & M ( \| A_1 x_1^k - A_1 x_1^{k+1} \|^2 + \| A_2 x_2^k - A_2 x_2^{k+1} \|^2 + \| x_3^k - x_3^{k+1} \|^2 ),
\end{eqnarray}
where $M:= \min\left\{ \frac{\gamma}{2}, \frac{\gamma + \sigma}{2} - \frac{L^2}{\gamma} \right\}$. Since $\gamma$ satisfies \eqref{gamma-1}, we have $M > 0$.

Then we prove that $\LCal_\gamma(w^k)$ is uniformly lower bounded. Since $f_1$, $f_2$ and $f_3$ are all lower bounded, we have
\begin{eqnarray}\label{inequality-augmented-lower}
& & \LCal_\gamma\left(x_1^{k+1},x_2^{k+1}, x_3^{k+1},\lambda^{k+1}\right) \nonumber \\
& \geq & f_1(x_1^{k+1}) + f_2(x_2^{k+1}) + f_3\left(b - \sum\limits_{i=1}^2 A_i x_i^{k+1}\right) + \frac{\gamma-L}{2}\left\|\sum_{i=1}^2 A_i x_i^{k+1} + x_3^{k+1} -b\right\|^2 \nonumber \\
& > & f_1^* + f_2^* + f_3^* := L^*,
\end{eqnarray}
where the first inequality holds from the convexity of $f_3$ and the Lipschitz continuity of $\nabla f_3$. By combining \eqref{L-decrease-M} and \eqref{inequality-augmented-lower}, for any integer $K>0$ we have
\begin{eqnarray*}
& & \sum\limits_{k=0}^K \left( \left\| A_1 x_1^k - A_1 x_1^{k+1} \right\|^2 + \left\| A_2 x_2^k - A_2 x_2^{k+1} \right\|^2 + \left\| x_3^k - x_3^{k+1} \right\|^2 \right) \\
& \leq & \frac{1}{M}\sum\limits_{k=0}^K \left(\LCal_\gamma(x_1^k,x_2^k,x_3^k,\lambda^k) - \LCal_\gamma(x_1^{k+1},x_2^{k+1},x_3^{k+1},\lambda^{k+1}) \right)  \\
& \leq & \frac{1}{M} \left( \LCal_\gamma(x_1^0,x_2^0,x_3^0,\lambda^0) - L^*\right).
\end{eqnarray*}
Letting $K\rightarrow+\infty$ yields
\[\sum\limits_{k=0}^\infty \left( \left\| A_1 x_1^k - A_1 x_1^{k+1} \right\|^2 + \left\| A_2 x_2^k - A_2 x_2^{k+1} \right\|^2 + \left\| x_3^k - x_3^{k+1} \right\|^2 \right) < +\infty,\]
which combining with \eqref{lambda-bound-x} yields
\be\label{sequence-convergence}
\lim_{k\rightarrow\infty}\|A_1x_1^k-A_1x_1^{k+1}\|=0, \lim_{k\rightarrow\infty}\|A_2x_2^k-A_2x_2^{k+1}\|=0, \lim_{k\rightarrow\infty}\| x_3^k- x_3^{k+1}\|=0, \lim_{k\rightarrow\infty}\|\lambda^k-\lambda^{k+1}\|=0.
\ee

Since $\LCal_\gamma(x_1^k,x_2^k,x_3^k;\lambda^k)$ is non-increasing and lower bounded, it follows that $\LCal_\gamma(x_1^k,x_2^k,x_3^k;\lambda^k)$ is convergent.
Finally, we prove that $\{\left(x_1^k,x_2^k,x_3^k,\lambda^k\right)\}$ is a bounded sequence. Note that \eqref{inequality-augmented-lower} and the coerciveness of $f_1+\textbf{1}_{\XCal_1}$ and $f_2+\textbf{1}_{\XCal_2}$ imply that
$\left\{\left(x_1^k,x_2^k\right): k=0,1,2,\ldots\right\}$ is a bounded sequence. This together with the updating formula of $\lambda^{k+1}$ and \eqref{sequence-convergence} yields the boundedness of $x_3^k$. Moreover, this combining with \eqref{opt-x-3-lambda} gives the boundedness of $\lambda^k$. Hence, $\left\{\left(x_1^k,x_2^k,x_3^k,\lambda^k\right): k=0,1,2,\ldots\right\}$ is a bounded sequence.

Therefore, there exists a limit point $(x_1^*,x_2^*,x_3^*,\lambda^*)$ and a subsequence $\{k_q\}$ such that
\begin{displaymath}
\lim\limits_{q\rightarrow\infty}x_i^{k_q} = x_i^*, i=1,2,3, \quad \lim\limits_{q\rightarrow\infty}\lambda^{k_q} = \lambda^*.
\end{displaymath}
From \eqref{sequence-convergence} we know
\begin{displaymath}
\lim\limits_{q\rightarrow\infty}x_i^{k_q+1} = x_i^*, i=1,2,3, \quad \lim\limits_{q\rightarrow\infty}\lambda^{k_q+1} = \lambda^*.
\end{displaymath}
Since $\LCal_\gamma(x_1^k,x_2^k,x_3^k;\lambda^k)$ is convergent, we know that
\be\label{L-convergent}\lim_{k\to\infty}\LCal_\gamma(x_1^k,x_2^k,x_3^k;\lambda^k) =
\LCal_\gamma(x_1^*,x_2^*,x_3^*;\lambda^*).
\ee
By combining the update of $x_3$ and $\lambda$, \eqref{opt-x-1-lambda} and \eqref{opt-x-2-lambda},
we know the following relations for any $x_1\in\XCal_1$ and $x_2\in\XCal_2$:
\begin{eqnarray*}
f_1(x_1) - f_1(x_1^{k_q+1}) + \left(x_1 - x_1^{k_q+1}\right)^\top\left[ - A_1^\top\lambda^{k_q+1} +\gamma A_1^\top\left( A_2(x_2^{k_q} - x_2^{k_q+1})+ (x_3^{k_q} - x_3^{k_q+1})\right)\right] & \geq & 0, \\
f_2(x_2) - f_2(x_2^{k_q+1}) + \left(x_2 - x_2^{k_q+1}\right)^\top\left[ - A_2^\top\lambda^{k_q+1} +\gamma A_2^\top\left( x_3^{k_q} - x_3^{k_q+1}\right)\right] & \geq & 0, \\
\nabla f_3(x_3^{k_q+1}) - \lambda^{k_q+1} & = & 0, \\
A_1 x_1^{k_q+1} + A_2 x_2^{k_q+1} + x_3^{k_q+1} - b - \frac{1}{\gamma}\left(\lambda^{k_q} - \lambda^{k_q+1}\right) & = & 0.
\end{eqnarray*}
Letting $q\rightarrow +\infty$, and using \eqref{sequence-convergence} and the lower semi-continuity of $f_1$ and $f_2$, we have
the following relations for any $x_1\in\XCal_1$ and $x_2\in\XCal_2$:
\begin{eqnarray*}
f_1(x_1) - f_1(x_1^*) - \left(x_1 - x_1^*\right)^\top ( A_1^\top\lambda^*  ) & \geq & 0, \\
f_2(x_2) - f_2(x_2^*) - \left(x_2 - x_2^*\right)^\top ( A_2^\top\lambda^*  ) & \geq & 0, \\
\nabla f_3(x_3^*) - \lambda^* & = & 0, \\
A_1 x_1^* + A_2 x_2^* + x_3^* - b & = & 0.
\end{eqnarray*}
Therefore, $\left(x_1^*,x_2^*,x_3^*,\lambda^*\right)$ satisfies the optimality conditions of problem \eqref{prob:P-uncons} and is an optimal primal and dual pair of problem \eqref{prob:P-uncons}.
Moreover, we have
\begin{eqnarray*}
& & \| A_1 x_1^k + A_2 x_2^k + x_3^k - b  \|  = \frac{1}{\gamma} \| \lambda^{k-1} - \lambda^k \| \rightarrow 0, \quad\mbox{ when } k\rightarrow \infty,
\end{eqnarray*}
and
\begin{eqnarray*}
& & \left| f(x_1^k) + f_2(x_2^k) + f_3(x_3^k) - f^*\right|  \nonumber \\
& \leq & \left| \LCal_\gamma ( x_1^k,x_2^k,x_3^k,\lambda^k ) - \LCal_\gamma ( x_1^*,x_2^*,x_3^*,\lambda^* ) \right| + \|\lambda^k \| \cdot\| A_1 x_1^k + A_2 x_2^k + x_3^k -b \|  \nonumber \\
& & + \frac{\gamma}{2} \| A_1 x_1^k + A_2 x_2^k + x_3^k -b \|^2\rightarrow 0, \quad \mbox{ when } k\to \infty,
\end{eqnarray*}
where we used \eqref{L-convergent}. Thus, we proved \eqref{thm-obj-feasibility}.

If Assumption \ref{assumption-4} holds, we can prove that the whole sequence of $(x_1^k,x_2^k,x_3^k;\lambda^k)$ converges to an optimal primal and dual pair of problem \eqref{prob:P-uncons}. Specifically, we have (note that in Case (i) of Assumption \ref{assumption-4}, $\sigma_1$ and $\sigma_2$ can be 0):
\begin{eqnarray}
& & \LCal_\gamma\left(x_1^k, x_2^k, x_3^k, \lambda^k\right) - \LCal_\gamma\left(x_1^*, x_2^*, x_3^*, \lambda^*\right)\nonumber\\
& \geq & \sum_{i=1}^2 \left(x_i^k - x_i^*\right)^\top \left(g_i^* - A_i^\top\lambda^*\right) + \frac{\sigma_1}{2}\left\| x_1^k - x_1^*\right\|^2 + \frac{\sigma_2}{2}\left\| x_2^k - x_2^*\right\|^2 + \frac{\sigma}{2}\left\| x_3^k - x_3^*\right\|^2  \nonumber\\
& & - \left\langle \lambda^k - \lambda^*, A_1 x_1^k + A_2 x_2^k + x_3^k - b\right\rangle + \frac{\gamma}{2}\left\| A_1 x_1^k + A_2 x_2^k + x_3^k - b\right\|^2 \nonumber\\
& \geq & \frac{\sigma_1}{2}\left\| x_1^k - x_1^*\right\|^2 + \frac{\sigma_2}{2}\left\| x_2^k - x_2^*\right\|^2 + \frac{\sigma}{2}\left\| x_3^k - x_3^*\right\|^2 - \frac{1}{\gamma}\left\langle \lambda^k - \lambda^*, \lambda^k - \lambda^{k+1}\right\rangle + \frac{1}{2\gamma}\left\| \lambda^k - \lambda^{k+1}\right\|^2 \nonumber\\
& \geq & \frac{\sigma_1}{2}\left\| x_1^k - x_1^*\right\|^2 + \frac{\sigma_2}{2}\left\| x_2^k - x_2^*\right\|^2 + \frac{\sigma}{2}\left\| x_3^k - x_3^*\right\|^2 - \frac{1}{2\gamma}\left\| \lambda^k - \lambda^*\right\|^2 \nonumber\\
& \geq & \frac{\sigma_1}{2}\left\| x_1^k - x_1^*\right\|^2 + \frac{\sigma_2}{2}\left\| x_2^k - x_2^*\right\|^2 + \left(\frac{\sigma}{2}-\frac{L^2}{2\gamma}\right)\left\| x_3^k - x_3^*\right\|^2.\label{L-decrease-new}
\end{eqnarray}
Since $\gamma$ satisfies \eqref{gamma-1}, we know $\frac{\sigma}{2}-\frac{L^2}{2\gamma} > 0$. Therefore, we have $x_3^k\rightarrow x_3^*$, which further implies $\lambda^k\rightarrow \lambda^*$. In the Case (i) of Assumption \ref{assumption-4}, since $A_1(x_1^k-x_1^*)+A_2(x_2^k-x_2^k)\rightarrow 0$ and $[A_1 \ A_2]$ is of full column rank, we know $x_1^k\rightarrow x_1^*$ and $x_2^k\rightarrow x_2^*$. In the Case (ii) of Assumption \ref{assumption-4}, since $\sigma_1>0$ and $\sigma_2>0$, \eqref{L-decrease-new} directly implies $(x_1^k,x_2^k,x_3^k)\rightarrow (x_1^*,x_2^*,x_3^*)$.

Now we prove part 2. We first assume that $\gamma\in\left(\sqrt{\sigma^2 + \frac{2L^2}{\eta_1 - 2}} - \sigma, \frac{4\sigma}{\eta_1}\right]$ for some $\eta_1>2$ such that $\sqrt{\sigma^2 + \frac{2L^2}{\eta_1 - 2}} - \sigma < \frac{4\sigma}{\eta_1}$. For any $(x_1^*,x_2^*,x_3^*,\lambda^*)\in\Omega^*$, combining \eqref{opt-x-1-lambda}-\eqref{opt-x-3-lambda} with \eqref{kkt} yields
\begin{eqnarray*}
& & \frac{1}{\gamma}\left(\lambda^k - \lambda^{k+1}\right)^\top\left(\lambda^{k+1}-\lambda^*\right) - \gamma \left( A_1 x_1^{k+1} - A_1 x_1^*\right)^\top \left( (A_2 x_2^k - A_2 x_2^{k+1}) + (x_3^k - x_3^{k+1}) \right) \\
& & - \gamma \left( A_2 x_2^{k+1} - A_2 x_2^*\right)^\top \left( x_3^k - x_3^{k+1} \right) \\
& \geq & \sigma \left\| x_3^{k+1} - x_3^* \right\|^2,
\end{eqnarray*}
which can be reduced to
\begin{eqnarray*}
& & \frac{1}{\gamma}\left( \lambda^k - \lambda^{k+1} \right)^\top \left(\lambda^{k+1}-\lambda^*\right) - \left(\lambda^k - \lambda^{k+1}\right)^\top \left( ( A_2 x_2^k - A_2 x_2^{k+1}) + ( x_3^k  - x_3^{k+1} ) \right) \\
& & + \gamma \left( A_2 x_2^{k+1} - A_2 x_2^*\right)^\top \left( A_2 x_2^{k} - A_2 x_2^{k+1} \right) + \gamma \left( x_3^{k+1} - x_3^*\right)^\top \left( (A_2 x_2^k - A_2 x_2^{k+1}) + (x_3^k - x_3^{k+1})\right) \\
& \geq & \sigma\left\| x_3^{k+1} - x_3^*\right\|^2.
\end{eqnarray*}
Combining this with \eqref{inequality-x-3-primal-dual} yields
\begin{eqnarray}\label{inequality-norm-decrease-1}
& & \frac{1}{\gamma}\left(\lambda^k - \lambda^{k+1}\right)^\top\left(\lambda^{k+1}-\lambda^*\right) + \gamma \left(A_2 x_2^{k+1} - A_2 x_2^*\right)^\top \left( A_2 x_2^{k} - A_2 x_2^{k+1} \right) \nonumber \\
& & + \gamma \left( x_3^{k+1} - x_3^*\right)^\top \left( x_3^k - x_3^{k+1} \right) \nonumber \\
& \geq & \sigma \left\| x_3^{k+1} - x_3^* \right\|^2 + \sigma \left\| x_3^k - x_3^{k+1} \right\|^2 + \left(\lambda^k - \lambda^{k+1}\right)^\top\left( A_2 x_2^{k} - A_2 x_2^{k+1} \right) \nonumber \\
& & - \gamma \left( x_3^{k+1} - x_3^* \right)^\top \left( A_2 x_2^{k} - A_2 x_2^{k+1} \right).
\end{eqnarray}
Now by applying \eqref{identity-3} to the three terms on the left hand side of \eqref{inequality-norm-decrease-1} we get
\begin{eqnarray}\label{inequality-norm-decrease-2}
& & \left[\frac{1}{2\gamma}\left\| \lambda^k - \lambda^* \right\|^2 + \frac{\gamma}{2} \left\| A_2 x_2^k - A_2 x_2^* \right\|^2 + \frac{\gamma}{2}\left\| x_3^k - x_3^*\right\|^2 \right] \nonumber \\
& & - \left[\frac{1}{2\gamma} \left\| \lambda^{k+1} - \lambda^* \right\|^2 + \frac{\gamma}{2} \left\| A_2 x_2^{k+1} - A_2 x_2^* \right\|^2 + \frac{\gamma}{2} \left\| x_3^{k+1} - x_3^* \right\|^2 \right] \nonumber \\
& \geq & \sigma\left\| x_3^{k+1} - x_3^* \right\|^2 + \sigma\left\| x_3^{k+1} - x_3^k\right\|^2 + \frac{1}{2\gamma} \left\|\lambda^{k+1} - \lambda^k \right\|^2 + \frac{\gamma}{2} \left\| A_2 x_2^{k+1} - A_2 x_2^k \right\|^2 +  \frac{\gamma}{2} \left\| x_3^{k+1} - x_3^k \right\|^2 \nonumber \\
& & + \left(\lambda^k - \lambda^{k+1}\right)^\top\left( A_2 x_2^{k} - A_2 x_2^{k+1} \right) - \gamma \left( x_3^{k+1} - x_3^*\right)^\top \left( A_2 x_2^k - A_2 x_2^{k+1} \right).
\end{eqnarray}
For any given $\eta_1>2$, we have
\begin{eqnarray}\label{inequality-couple-1}
-\gamma\left( x_3^{k+1} - x_3^* \right)^\top \left( A_2 x_2^k - A_2 x_2^{k+1} \right) \geq -\frac{\gamma\eta_1}{4}\left\| x_3^{k+1} - x_3^* \right\|^2 - \frac{\gamma}{\eta_1}\left\| A_2 x_2^{k+1} - A_2 x_2^k \right\|^2 ,
\end{eqnarray}
and
\begin{eqnarray}\label{inequality-couple-2}
& & \frac{\eta_1}{2\gamma(\eta_1-2)} \left\|\lambda^{k+1} - \lambda^k \right\|^2 + \left(\lambda^{k+1} - \lambda^k\right)^\top \left( A_2 x_2^{k+1} - A_2 x_2^k \right) + \frac{\gamma(\eta_1-2)}{2\eta_1} \left\| A_2 x_2^{k+1} - A_2 x_2^k \right\|^2
\nonumber \\
& = & \left\|\sqrt{\frac{\eta_1}{2\gamma(\eta_1-2)}}\left( \lambda^{k+1} - \lambda^k \right) + \sqrt{\frac{\gamma(\eta_1-2)}{2\eta_1}}\left( A_2 x_2^{k+1} - A_2 x_2^k \right)\right\|^2.
\end{eqnarray}
By combining \eqref{lambda-bound-x}, \eqref{inequality-couple-1}, \eqref{inequality-couple-2} and \eqref{inequality-norm-decrease-2}, we get
\begin{eqnarray}
& & \left[\frac{1}{2\gamma} \left\| \lambda^k - \lambda^* \right\|^2 + \frac{\gamma}{2} \left\| A_2 x_2^k - A_2 x_2^* \right\|^2 + \frac{\gamma}{2} \left\| x_3^k - x_3^* \right\|^2 \right]
\nonumber \\
& & - \left[\frac{1}{2\gamma} \left\| \lambda^{k+1} - \lambda^* \right\|^2 + \frac{\gamma}{2} \left\| A_2 x_2^{k+1} - A_2 x_2^* \right\|^2 + \frac{\gamma}{2} \left\| x_3^{k+1} - x_3^* \right\|^2 \right] \nonumber \\
& \geq & \left(\sigma - \frac{\gamma\eta_1}{4}\right)\left\| x_3^{k+1} - x_3^* \right\|^2 + \left(\sigma + \frac{\gamma}{2} - \frac{L^2}{\gamma(\eta_1-2)} \right) \left\| x_3^{k+1} - x_3^k \right\|^2 \nonumber \\
& & + \left\|\sqrt{\frac{\eta_1}{2\gamma(\eta_1-2)}}\left( \lambda^{k+1} - \lambda^k \right) + \sqrt{\frac{\gamma(\eta_1-2)}{2\eta_1}}\left( A_2 x_2^{k+1} - A_2 x_2^k \right)\right\|^2 \nonumber \\
& \geq & 0,\label{telescopic-gamma-eta1}
\end{eqnarray}
where the second inequality holds because $\gamma\in\left(\sqrt{\sigma^2 + \frac{2L^2}{\eta_1 - 2}} - \sigma, \frac{4\sigma}{\eta_1}\right]$ implies that $$\sigma - \frac{\eta_1\gamma}{4}\geq 0, \quad \sigma + \frac{\gamma}{2} - \frac{L^2}{\gamma(\eta_1-2)}>0.$$

Furthermore, \eqref{telescopic-gamma-eta1} implies $\| x_3^{k+1} - x_3^k \|\rightarrow 0$ and hence $\| \lambda^{k+1} - \lambda^k\|\rightarrow 0$ because of $\nabla f_3(x_3^k) = \lambda^k$, and $\| A_2 x_2^{k+1} - A_2 x_2^k\|\rightarrow 0$ since
\[
\left\|\sqrt{\frac{\eta_1}{2\gamma(\eta_1-2)}}\left( \lambda^{k+1} - \lambda^k \right) + \sqrt{\frac{\gamma(\eta_1-2)}{2\eta_1}}\left( A_2 x_2^{k+1} - A_2 x_2^k \right)\right\|\rightarrow 0.
\]
Moreover, the sequence $\frac{1}{2\gamma} \left\| \lambda^k - \lambda^* \right\|^2 + \frac{\gamma}{2} \left\| A_2 x_2^k - A_2 x_2^* \right\|^2 + \frac{\gamma}{2} \left\| x_3^k - x_3^* \right\|^2$ is non-increasing, and this implies that $\left\{\left(A_2 x_2^k,x_3^k,\lambda^k\right): k=0,1,2,\ldots\right\}$ is bounded. Since $A_1$ and $A_2$ both have full column rank, we conclude that $\left\{\left(x_1^k, x_2^k,x_3^k,\lambda^k\right): k=0,1,2,\ldots\right\}$ is a bounded sequence.

Therefore, there exists a limit point $\left(\bar {x}_1,\bar{x}_2,\bar{x}_3,\bar{\lambda}\right)$ and a subsequence $\{k_q\}$ such that
\begin{displaymath}
\lim\limits_{q\rightarrow\infty}x_i^{k_q} = \bar{x}_i, i=1,2,3, \quad \lim\limits_{q\rightarrow\infty}\lambda^{k_q} = \bar{\lambda}.
\end{displaymath}
By $\| A_2 x_2^{k+1} - A_2 x_2^k \|\rightarrow 0$, $\| x_3^{k+1} - x_3^k \|\rightarrow 0$ and $\| \lambda^{k+1} - \lambda^k\|\rightarrow 0$, we have
\begin{displaymath}
\lim\limits_{q\rightarrow\infty}x_i^{k_q+1} = \bar{x}_i, i=1,2,3, \quad \lim\limits_{q\rightarrow\infty}\lambda^{k_q+1} = \bar{\lambda}.
\end{displaymath}
By the same argument as in the above case, we conclude that $\left(\bar{x}_1,\bar{x}_2,\bar{x}_3,\bar{\lambda}\right)$ is an optimal primal and dual pair of \eqref{prob:P-uncons}.

Finally, we prove that the whole sequence $(x_1^k, x_2^k,x_3^k,\lambda^k)$ converges to $(\bar{x}_1, \bar{x}_2, \bar{x}_3, \bar{\lambda})$. It suffices to prove that $(A_1 x_1^k, A_2 x_2^k,x_3^k,\lambda^k)$ converges to $(A_1 \bar{x}_1, A_2 \bar{x}_2, \bar{x}_3, \bar{\lambda})$ since $A_1$ and $A_2$ both have full column rank.
Note that since $(\bar{x}_1, \bar{x}_2, \bar{x}_3, \bar{\lambda})$ is an optimal primal and dual pair of \eqref{prob:P-uncons}, \eqref{telescopic-gamma-eta1} holds with $(x_2^*,x_3^*,\lambda^*)$ replaced by $(\bar{x}_2,\bar{x}_3,\bar{\lambda})$. Therefore,
$\frac{1}{2\gamma} \| \lambda^{k} - \bar{\lambda}\|^2 + \frac{\gamma}{2} \| A_2 x_2^{k} - A_2 \bar{x}_2\|^2 + \frac{\gamma}{2} \| x_3^{k} - \bar{x}_3 \|^2$ is non-increasing. Moreover, we have
$\frac{1}{2\gamma} \| \lambda^{k_q} - \bar{\lambda} \|^2 + \frac{\gamma}{2} \| A_2 x_2^{k_q} - A_2 \bar{x}_2\|^2 + \frac{\gamma}{2} \| x_3^{k_q} - \bar{x}_3\|^2 \rightarrow 0$. Therefore, it follows that
\begin{displaymath}
\frac{1}{2\gamma} \| \lambda^k - \bar{\lambda}  \|^2 + \frac{\gamma}{2} \| A_2 x_2^k - A_2 \bar{x}_2\|^2 + \frac{\gamma}{2} \| x_3^k - \bar{x}_3 \|^2 \rightarrow 0,
\end{displaymath}
i.e., the whole sequence of $(A_2x_2^k,x_3^k,\lambda^k)$ converges to $(A_2\bar{x}_2,\bar{x}_3,\bar{\lambda})$.
Furthermore, $\| A_1 x_1^k - A_1 \bar{x}_1\|\rightarrow 0$ by using the update formula of $\lambda^{k+1}$.

Now we assume $\gamma\in\left(0, \min\left\{\frac{4\sigma}{\eta_2}, \frac{\sigma(\eta_2 - 2)}{4\eta_2} + \sqrt{\frac{\sigma^2(\eta_2 - 2)^2}{16\eta_2^2} + \frac{\sigma^2(\eta_2 - 2)}{4\eta_2}}\right\}\right)$ for arbitrarily chosen $\eta_2>2$.
Using similar arguments as in the case $\gamma\in\left(\sqrt{\sigma^2 + \frac{2L^2}{\eta_1 - 2}} - \sigma, \frac{4\sigma}{\eta_1}\right]$, the following inequalities hold for any given $\eta_2>2$ and $\epsilon > \frac{2\eta_2}{\eta_2-2}$:
\begin{equation}\label{inequality-couple-3}
-\gamma\left( x_3^{k+1} - x_3^* \right)^\top \left( A_2 x_2^k - A_2 x_2^{k+1} \right)\geq -\frac{\gamma\eta_2}{4} \left\| x_3^{k+1} -  x_3^* \right\|^2 - \frac{\gamma}{\eta_2} \left\| A_2 x_2^{k} - A_2 x_2^{k+1} \right\|^2,
\end{equation}
and
\begin{eqnarray}\label{inequality-couple-4}
& & \left( \lambda^k - \lambda^{k+1} \right)\left(A_2 x_2^k - A_2 x_2^{k+1}\right)^{\top} \nonumber \\
& \geq & -\frac{\gamma}{\epsilon}\left\| A_2 x_2^{k+1} - A_2 x_2^k \right\|^2 - \frac{\gamma\epsilon}{2}\left\| x_3^k - x_3^{k-1}\right\|^2 - \frac{\gamma\epsilon}{2}\left\| x_3^{k+1} - x_3^k\right\|^2.
\end{eqnarray}
It follows from \eqref{strong_convex_2} and \eqref{opt-x-3-lambda} that
\begin{equation}\label{inequality-couple-5}
\left\| \lambda^{k+1} - \lambda^* \right\| \geq \sigma\left\| x_3^{k+1} - x_3^* \right\|.
\end{equation}

Therefore, we conclude from \eqref{inequality-couple-3}-\eqref{inequality-couple-5} and \eqref{inequality-norm-decrease-2} that
\begin{eqnarray*}
& & \left[ \frac{1}{2\gamma} \left\| \lambda^k - \lambda^* \right\|^2 + \frac{\gamma}{2} \left\| A_2 x_2^k - A_2 x_2^* \right\|^2 + \frac{\gamma}{2} \left\| x_3^k - x_3^* \right\|^2 + \frac{\gamma\epsilon}{2}\left\| x_3^k - x_3^{k-1}\right\|^2\right] \\
& & - \left[\frac{1}{2\gamma} \left\| \lambda^{k+1}- \lambda^* \right\|^2 + \frac{\gamma}{2} \left\| A_2 x_2^{k+1} - A_2 x_2^* \right\|^2 + \frac{\gamma}{2} \left\| x_3^{k+1} - x_3^* \right\|^2 + \frac{\gamma\epsilon}{2}\left\| x_3^{k+1} - x_3^k\right\|^2\right] \\
& \geq & \left(\sigma - \frac{\eta_2\gamma}{4}\right) \left\| x_3^{k+1} - x_3^* \right\|^2 + \left(\sigma + \frac{\sigma^2}{2\gamma} - \gamma\epsilon\right) \left\| x_3^{k+1} - x_3^k \right\|^2 + \left(\frac{\gamma}{2} - \frac{\gamma}{\eta_2} - \frac{\gamma}{\epsilon}\right) \left\| A_2 x_2^{k+1} - A_2 x_2^k \right\|^2 \\
& \geq & \left(\sigma + \frac{\sigma^2}{2\gamma} - \gamma\epsilon\right) \left\| x_3^{k+1} - x_3^k \right\|^2 + \left(\frac{\gamma}{2} - \frac{\gamma}{\eta_2} - \frac{\gamma}{\epsilon}\right) \left\| A_2 x_2^{k+1} - A_2 x_2^k \right\|^2 \nonumber \\
& \geq & 0,
\end{eqnarray*}
where the second and third inequalities hold because $\gamma\in\left(0, \min\left\{\frac{4\sigma}{\eta_2}, \frac{\sigma(\eta_2 - 2)}{4\eta_2} + \sqrt{\frac{\sigma^2(\eta_2 - 2)^2}{16\eta_2^2} + \frac{\sigma^2(\eta_2 - 2)}{4\eta_2}}\right\}\right)$ for any $\eta_2>2$ implies
$$0 < \gamma \leq \frac{4\sigma}{\eta_2}, \quad \frac{\gamma}{2} - \frac{\gamma}{\eta_2} - \frac{\gamma}{\epsilon} > 0, \quad \sigma + \frac{\sigma^2}{2\gamma} - \gamma\epsilon>0.$$

This implies $\| x_3^{k+1} - x_3^k \|\rightarrow 0$, $\| A_2 x_2^{k+1} - A_2 x_2^k \|\rightarrow 0$, and hence $\| \lambda^{k+1} - \lambda^k\|\rightarrow 0$. This also implies the sequence $\frac{1}{2\gamma} \left\| \lambda^k - \lambda^* \right\|^2 + \frac{\gamma}{2} \left\| A_2 x_2^k - A_2 x_2^* \right\|^2 + \frac{\gamma}{2} \left\| x_3^k - x_3^* \right\|^2 + \frac{\gamma\epsilon}{2}\left\| x_3^k - x_3^{k-1} \right\|^2$ is non-increasing, which further implies that $\left\{\left(A_2 x_2^k,x_3^k,\lambda^k\right): k=0,1,2,\ldots\right\}$ is bounded.  Since $A_1$ and $A_2$ both have full column rank, we conclude that $\left\{\left(x_1^k, x_2^k,x_3^k,\lambda^k\right): k=0,1,2,\ldots\right\}$ is a bounded sequence.

Finally, using similar arguments as in the case $\gamma\in\left(\sqrt{\sigma^2 + \frac{2L^2}{\eta_1 - 2}} - \sigma, \frac{4\sigma}{\eta_1}\right]$ it is easy to prove that the whole sequence of $\left\{\left(x_1^k, x_2^k,x_3^k,\lambda^k\right): k=0,1,2,\ldots\right\}$ converges to $(x_1^*,x_2^*,x_3^*,\lambda^*)$. We omit the details here for succinctness.
\end{proof}

\begin{remark}
We remark here that there exist works that show the whole sequence convergence of 2-block ADMM for even nonconvex problems, but they usually require some other assumptions such as the Kurdyka-{\L}ojasiewicz property (see, e.g., \cite{Li-Pong-2015}).
\end{remark}

The following proposition shows that the interval in \eqref{assumption} equals $(0,\infty)$ if the condition number of $f_3$ is in $[1,1.0798)$.
\begin{proposition}\label{prop:gamma-free}
If the condition number of $f_3$, i.e., $\kappa := L/\sigma$, is in $[1,1.0798)$, then there exist $\eta_1, \eta_2\in (2,\infty)$, such that \eqref{assumption} reduces to $\gamma\in(0,\infty)$. That is, $\gamma$ can be freely chosen in $(0,\infty)$.
\end{proposition}
\begin{proof}
Without loss of generality, we can assume that $\sigma=1$. Therefore, $L=\kappa$. By letting
$\frac{4}{\eta_2} = \frac{\eta_2-2}{4\eta_2}+\sqrt{\frac{(\eta_2-2)^2}{16\eta_2^2}+\frac{\eta_2-2}{4\eta_2}}$, we have $\eta_2=\sqrt{89}-3$. By letting $4/\eta_1=\kappa^2$, we have $\kappa <\sqrt{2}$ because of $\eta_1>2$.
In addition, we need $\sqrt{1+\frac{2\kappa^2}{\eta_1-2}}-1 < 4/\eta_2$, and we found $\kappa^2 < 1.1659$, i.e., $\kappa < 1.0798$ suffices.
\end{proof}

\begin{remark}
Note that $\kappa=1$ implies that $f_3(\cdot)=\half\|\cdot\|^2$. Therefore, 3-block ADMM globally converges for any $\gamma>0$ when it is applied to solve the RLSD problem.
\end{remark}



\section{Numerical Experiments}\label{sec:num}

While comparing 3-block ADMM with other methods is not the main focus of this paper, we shall present some numerical results in this section to gain some insights on the performance of these methods. 

{
\subsection{Two alternative ways for solving RLSD \eqref{sharing}}
Here we discuss two alternative approaches for solving RLSD \eqref{sharing} and then compare them with the 3-block ADMM. 
One natural way to solve \eqref{sharing} is to apply the block coordinate descent (BCD) method, where the iterates are updated as
\be\label{bcd}
\left\{\ba{lll} x_1^{k+1} & := & \argmin_{x_1\in\XCal_1} \ f_1(x_1) + \half\|A_1x_1+A_2x_2^k-b\|^2 \\
                x_2^{k+1} & := & \argmin_{x_2\in\XCal_2} \ f_2(x_2) + \half\|A_1x_1^{k+1}+A_2x_2-b\|^2. \ea\right.
\ee
The other way for solving \eqref{prob:P-uncons} is to apply the 2-block ADMM. Specifically, by grouping $(x_2,x_3)$ as one block variable, \eqref{prob:P-uncons} can be solved by 2-block ADMM as follows:
\be\label{admm-2-referee}
\left\{\ba{lll} x_1^{k+1} & := & \argmin_{x_1\in\XCal_1} \ \LCal_\gamma(x_1,x_2^k,x_3^k;\lambda^k) \\
               (x_2^{k+1},x_3^{k+1}) & := & \argmin_{x_2\in\XCal_2, x_3} \ \LCal_\gamma(x_1^{k+1},x_2,x_3;\lambda^k) \\
               \lambda^{k+1} & := & \lambda^k - \gamma(A_1x_1+A_2x_2+x_3-b). \ea\right.
\ee
Due to the special structure of \eqref{prob:P-uncons}, i.e., $f_3(x_3) = \half\|x_3\|^2$, the second subproblem in \eqref{admm-2-referee} is equivalent to
\be\label{admm-2-referee-sub}
\left\{\ba{lll} x_2^{k+1} & := & \argmin_{x_2\in\XCal_2} f_2(x_2) + \frac{\gamma}{2(1+\gamma)}\|A_1x_1^{k+1}+A_2x_2-b-\lambda^k/\gamma\|^2 \\
                x_3^{k+1} & := & \frac{1}{1+\gamma}(\lambda^k-\gamma(A_1x_1^{k+1}+A_2x_2^{k+1}-b)).\ea\right.
\ee
It is thus noted that both BCD \eqref{bcd} and the 2-block ADMM \eqref{admm-2-referee} have the same per-iteration complexity as 3-block ADMM for solving \eqref{sharing} and \eqref{prob:P-uncons}. Moreover, BCD \eqref{bcd} does not need any parameter, and 2-block ADMM \eqref{admm-2-referee} globally converges for any $\gamma>0$. As a result, both \eqref{bcd} and \eqref{admm-2-referee} are natural choices for solving \eqref{sharing} and \eqref{prob:P-uncons}.
We shall conduct some numerical comparisons of 3-block ADMM \eqref{admm-3}, 2-block ADMM \eqref{admm-2-referee} and BCD \eqref{bcd} for solving the stable principal component pursuit problem \eqref{spcp}. 
}

\subsection{Comparison of 3-Block ADMM with BCD and a 2-Block ADMM}

In this subsection, we report some numerical results on solving the SPCP problem \eqref{spcp}. It is noted that \eqref{spcp} can be solved by BCD, where the iterates are updated as
\be\label{spcp-bcd}
\left\{\ba{lll}
L^{k+1} & := & \argmin_L \ \beta_1\|L\|_* + \half\|M-L-S^k\|_F^2 \\
S^{k+1} & := & \argmin_S \ \beta_2\|S\|_1 + \half\|M-L^{k+1}-S\|_F^2.
\ea\right.
\ee
By equivalently reformulating \eqref{spcp} to
\be\label{spcp-reformulate} \min_{L,S,Z} \ \beta_1\|L\|_* + \beta_2\|S\|_1 + \half\|Z\|_F^2, \ \st, L+S+Z=M, \ee
we can apply both the 2-block ADMM \eqref{admm-2-referee} (denoted as ADMM-2) and the 3-block ADMM (denoted as ADMM-3) for solving it. To compare the performance of BCD, ADMM-2 and ADMM-3 for solving \eqref{spcp}, we tested them on some randomly created problems. The problems were created in similar manner as \cite{Lin-Chen-Wu-Ma-ADM-RPCA-2009}. For simplicity, we set $m=n$ in all the tested problems. The matrix $M$ was generated in the following way. For given $n$ and $r<n$, we set the targeting rank-$r$ matrix $L^*= L_1*L_2^\top$, where $L_1$ and $L_2$ are $n\times r$ matrices whose entries are i.i.d. Gaussian random variables drawn from $\mathcal{N}(0,1)$. For given sparsity $s$, the support of the targeting sparse matrix $S^*$ was chosen uniformly at random, and the $s$ nonzero entries were i.i.d. Gaussian random variables drawn from $\mathcal{N}(0,1)$. The entries of the noise matrix $Z^*$ follows i.i.d. Gaussian $\mathcal{N}(0,1)\times 10^{-8}$. Finally, we set $M=L^*+S^*+Z^*$. We set $\beta_1=0.005$ and $\beta_2 = \beta_1/\sqrt{n}$. We define the relative errors of $L$ and $S$ as
\[errL = \frac{\|L-L^*\|_F}{\|L^*\|_F}, \quad errS = \frac{\|S-S^*\|_F}{\|S^*\|_F},\]
and all three algorithms were terminated when $\max(errL,errS)<10^{-3}$, or the maximum number of iterations 20000 was reached. We tested the three algorithms for different $n$, $r$, $s$ and the results are reported in Tables \ref{tab:rpca-1} and \ref{tab:rpca-2}. Table \ref{tab:rpca-1} gives the results for $r=0.05n$ and $s=0.05n^2$ and Table \ref{tab:rpca-2} gives the results for $r=0.05n$ and $s=0.1n^2$. We set the initial Lagrange multiplier as 0. We tested two initial primal variables: $(L_0,S_0,Z_0) = (0,0,0)$ and $(L_0,S_0,Z_0) = (0,0,\gamma M/(1+\gamma))$. Note that the latter satisfies the second equation in \eqref{admm-2-referee-sub}. Based on the results in Tables \ref{tab:rpca-1} and \ref{tab:rpca-2}, we observed the following characteristics. First, the performance of BCD is very robust, but ADMM-3 with an appropriate $\gamma$ (i.e., $\gamma=0.7$) can outperform BCD in terms of number of iterations and CPU time required to reach the same error bounds $errL$ and $errS$, while ADMM-3 performs worse than BCD when $\gamma=1.2$. This indicates that the performance of ADMM-3 varies for different $\gamma$, and that if one knows how to choose the parameter $\gamma$, ADMM-3 can be faster than BCD, although the latter has no parameter to choose. Second, the performance of ADMM-2 seems to depend more on the initial primal variables than ADMM-3. For instance, when $(L_0,S_0,Z_0) = (0,0,0)$, ADMM-2 and ADMM-3 need almost the same number of iterations and CPU time to reach solutions with the same error margin; when $(L_0,S_0,Z_0) = (0,0,\gamma M/(1+\gamma))$, the performance of ADMM-3 is much better. This may seem counterintuitive at the first glance, because both ADMM-2 \eqref{admm-2-referee} and ADMM-3 \eqref{admm-3} globally converge for RLSD. We observe that updating $x_3$ in ADMM-3 \eqref{admm-3} requires the latest information of $x_1$ and $x_2$; while updating $x_3$ in ADMM-2 \eqref{admm-2-referee} only requires the latest information of $x_1$, because $x_2$ is also decided by $x_1$. It is our belief that this might explain why ADMM-3 is better than ADMM-2 in this case.

\begin{table}[t]\small
\centering
\vspace{-6em}
\begin{tabular}{|c|c|c|c|c|} \hline
$n$ & $errL$ & $errS$ & iter & CPU \\ \hline\hline


\multicolumn{5}{|c|} {$L_0=0$, $S_0=0$, $Z_0=0$} \\ \hline\hline

\multicolumn{5}{|c|} {BCD} \\ \hline

100 & 8.3450e-05 & 9.5182e-04 &  1380 & 6.97489 \\ \hline

200 & 6.5649e-05 & 9.5838e-04 &  1738 & 34.84547 \\ \hline

400 & 4.9232e-05 & 9.7616e-04 &  2175 & 188.19662 \\ \hline

\multicolumn{5}{|c|} {ADMM-3, $\gamma=0.7$} \\ \hline

100 & 8.3111e-05 & 9.4803e-04 &   966 & 4.89613 \\ \hline

200 & 6.2285e-05 & 9.0967e-04 &  1217 & 21.81870 \\ \hline

400 & 5.0408e-05 & 9.9912e-04 &  1522 & 132.73636 \\ \hline

\multicolumn{5}{|c|} {ADMM-2, $\gamma=0.7$} \\ \hline

100 & 8.6981e-05 & 9.9094e-04 &   966 & 4.64288 \\ \hline

200 & 6.4384e-05 & 9.4005e-04 &  1217 & 19.31321 \\ \hline

400 & 4.8454e-05 & 9.6096e-04 &  1523 & 130.00307 \\ \hline

\multicolumn{5}{|c|} {ADMM-3, $\gamma=1.2$} \\ \hline

100 & 8.3585e-05 & 9.5332e-04 &  1656 & 7.90678 \\ \hline

200 & 6.3907e-05 & 9.3315e-04 &  2086 & 34.85310 \\ \hline

400 & 4.9272e-05 & 9.7693e-04 &  2610 & 227.28830 \\ \hline

\multicolumn{5}{|c|} {ADMM-2, $\gamma=1.2$} \\ \hline

100 & 8.7394e-05 & 9.9552e-04 &  1656 & 8.07929 \\ \hline

200 & 6.5990e-05 & 9.6325e-04 &  2086 & 34.17833 \\ \hline

400 & 4.9858e-05 & 9.8837e-04 &  2610 & 203.37331 \\ \hline \hline

\multicolumn{5}{|c|} {$L_0=0$, $S_0=0$, $Z_0=\gamma M/(1+\gamma)$} \\ \hline\hline

\multicolumn{5}{|c|} {BCD} \\ \hline

100 & 1.0435e-04 & 9.7557e-04 &  1718 & 5.37146 \\ \hline

200 & 6.6891e-05 & 9.3881e-04 &  1754 & 24.71132 \\ \hline

400 & 4.9139e-05 & 9.7866e-04 &  2094 & 178.39425 \\ \hline

\multicolumn{5}{|c|} {ADMM-3, $\gamma=0.7$} \\ \hline

100 & 1.0503e-04 & 9.8249e-04 &  1002 & 3.42133 \\ \hline

200 & 6.5317e-05 & 9.1722e-04 &  1228 & 16.87290 \\ \hline

400 & 4.8614e-05 & 9.6818e-04 &  1466 & 113.56560 \\ \hline

\multicolumn{5}{|c|} {ADMM-2, $\gamma=0.7$} \\ \hline

100 & 1.0748e-04 & 9.2488e-04 &  5961 & 19.84691 \\ \hline

200 & 7.5803e-05 & 9.6968e-04 &  9974 & 152.25914 \\ \hline

400 & 5.0729e-05 & 9.7315e-04 & 18987 & 1643.33284 \\ \hline

\multicolumn{5}{|c|} {ADMM-3, $\gamma=1.2$} \\ \hline

100 & 1.0161e-04 & 9.4712e-04 &  1432 & 4.67480 \\ \hline

200 & 7.0248e-05 & 9.8482e-04 &  2104 & 33.59378 \\ \hline

400 & 5.0117e-05 & 9.9824e-04 &  2512 & 234.20045 \\ \hline

\multicolumn{5}{|c|} {ADMM-2, $\gamma=1.2$} \\ \hline

100 & 1.1461e-04 & 9.8808e-04 & 12870 & 42.19114 \\ \hline

200 & 6.6541e-03 & 9.0977e-02 & 20000 & 332.14747 \\ \hline

400 & 1.0979e-01 & 2.1819e+00 & 20000 & 1638.61065 \\ \hline

\end{tabular}\caption{The comparison results on SPCP with $r=0.05n$ and $s=0.05n^2$}\label{tab:rpca-1}
\end{table}

\begin{table}[t]\small
\centering
\vspace{-6em}
\begin{tabular}{|c|c|c|c|c|} \hline
$n$ & $errL$ & $errS$ & iter & CPU \\ \hline\hline


\multicolumn{5}{|c|} {$L_0=0$, $S_0=0$, $Z_0=0$} \\ \hline\hline

\multicolumn{5}{|c|} {BCD} \\ \hline

100 & 1.2020e-04 & 9.5076e-04 &  2254 & 9.22008 \\ \hline

200 & 9.3804e-05 & 9.6703e-04 &  2480 & 47.64641 \\ \hline

400 & 6.7383e-05 & 9.8883e-04 &  3191 & 309.17288 \\ \hline

\multicolumn{5}{|c|} {ADMM-3, $\gamma=0.7$} \\ \hline

100 & 1.1711e-04 & 9.2993e-04 &  1578 & 6.44356 \\ \hline

200 & 9.3552e-05 & 9.6463e-04 &  1736 & 33.24449 \\ \hline

400 & 6.6681e-05 & 9.7902e-04 &  2234 & 208.73834 \\ \hline

\multicolumn{5}{|c|} {ADMM-2, $\gamma=0.7$} \\ \hline

100 & 1.2174e-04 & 9.6117e-04 &  1578 & 6.69733 \\ \hline

200 & 9.5246e-05 & 9.8091e-04 &  1736 & 32.65144 \\ \hline

400 & 6.7352e-05 & 9.8840e-04 &  2234 & 226.20321 \\ \hline

\multicolumn{5}{|c|} {ADMM-3, $\gamma=1.2$} \\ \hline

100 & 1.2615e-04 & 9.9098e-04 &  2704 & 11.89260 \\ \hline

200 & 9.6866e-05 & 9.9645e-04 &  2975 & 56.88310 \\ \hline

400 & 6.7657e-05 & 9.9266e-04 &  3829 & 363.85171 \\ \hline

\multicolumn{5}{|c|} {ADMM-2, $\gamma=1.2$} \\ \hline

100 & 1.2351e-04 & 9.7312e-04 &  2705 & 11.98009 \\ \hline

200 & 9.5563e-05 & 9.8392e-04 &  2976 & 53.21315 \\ \hline

400 & 6.7136e-05 & 9.8538e-04 &  3830 & 327.37323 \\ \hline \hline

\multicolumn{5}{|c|} {$L_0=0$, $S_0=0$, $Z_0=\gamma M/(1+\gamma)$} \\ \hline\hline

\multicolumn{5}{|c|} {BCD} \\ \hline

100 & 1.4379e-04 & 9.8739e-04 &  2903 & 10.00203 \\ \hline

200 & 9.2676e-05 & 9.8801e-04 &  2571 & 37.43186 \\ \hline

400 & 6.8765e-05 & 9.8623e-04 &  3138 & 275.30181 \\ \hline

\multicolumn{5}{|c|} {ADMM-3, $\gamma=0.7$} \\ \hline

100 & 1.4392e-04 & 9.8819e-04 &  2032 & 6.90356 \\ \hline

200 & 9.0868e-05 & 9.7045e-04 &  1800 & 27.81380 \\ \hline

400 & 6.7811e-05 & 9.7323e-04 &  2197 & 214.79459 \\ \hline

\multicolumn{5}{|c|} {ADMM-2, $\gamma=0.7$} \\ \hline

100 & 1.4217e-04 & 9.7743e-04 &  6871 & 24.20091 \\ \hline

200 & 9.9873e-05 & 9.8780e-04 & 10867 & 197.24003 \\ \hline

400 & 7.1439e-05 & 9.8609e-04 & 16589 & 1487.18602 \\ \hline

\multicolumn{5}{|c|} {ADMM-3, $\gamma=1.2$} \\ \hline

100 & 1.4575e-04 & 9.9940e-04 &  3483 & 11.88450 \\ \hline

200 & 9.3381e-05 & 9.9486e-04 &  3085 & 56.96292 \\ \hline

400 & 6.9547e-05 & 9.9691e-04 &  3765 & 366.69591 \\ \hline

\multicolumn{5}{|c|} {ADMM-2, $\gamma=1.2$} \\ \hline

100 & 1.4801e-04 & 9.9357e-04 & 14135 & 49.91323 \\ \hline

200 & 7.7218e-03 & 7.7387e-02 & 20000 & 343.15011 \\ \hline

400 & 1.0250e-01 & 1.4256e+00 & 20000 & 1860.90799 \\ \hline

\end{tabular}\caption{The comparison results of SPCP with $r=0.05n$ and $s=0.1n^2$}\label{tab:rpca-2}
\end{table}

\section{Conclusions}

Motivated by the fact that the 2-block ADMM globally converges for any penalty parameter $\gamma>0$, we studied in this paper the global convergence of the 3-block ADMM. As there exists a counter-example showing that the 3-block ADMM can diverge if no further condition is imposed, it
is natural to look for
sufficient conditions that can guarantee the convergence of the 3-block ADMM. However, the existing results on sufficient conditions usually require $\gamma$ to be smaller than a certain bound, which is usually very small and therefore not practically efficient. In this paper, we showed that the 3-block ADMM globally converges for any $\gamma>0$ when if $A_3=I$, $x_3\in\mathcal{X}_3$ is absent, and more importantly, the condition number of $f_3$ is in $[1,1.0798)$; that is, the 3-block ADMM is parameter-unrestricted for this class of problems. 

%

\bibliographystyle{plain}
\bibliography{admm3-free-penalty}

\end{document}